# HAUSDORFF DIMENSIONS FOR SLE$_6$

By Vincent Beffara

*Université Paris-Sud and Institut Mittag-Leffler*

We prove that the Hausdorff dimension of the trace of SLE$_6$ is almost surely 7/4 and give a more direct derivation of the result (due to Lawler–Schramm–Werner) that the dimension of its boundary is 4/3. We also prove that, for all $\kappa < 8$, the SLE$_\kappa$ trace has cut-points.

**Contents**

**0. Introduction.** It has been conjectured by theoretical physicists that various lattice models in statistical physics (such as percolation, Potts model,

Ising, uniform spanning trees), taken at their critical point, have a continuous conformally invariant scaling limit when the mesh of the lattice tends to 0. Recently, Schramm [16] introduced a family of random processes he called Stochastic Loewner Evolutions (or SLE), which are the only possible conformally invariant scaling limits of random cluster interfaces (which are very closely related to all above-mentioned models).

An SLE process is defined using the usual Loewner equation, where the driving function is a time-changed Brownian motion. More specifically, in the present paper we will be mainly concerned with SLE in the upper-half plane (sometimes called chordal SLE), defined by the following PDE:

$$\partial_t g_t(z) = \frac{2}{g_t(z) - \sqrt{\kappa} B_t}, \qquad g_0(z) = z, \tag{1}$$

where $(B_t)$ is a standard Brownian motion on the real line and $\kappa$ is a positive parameter. It can be shown that this equation defines a family $(g_t)$ of conformal mappings from simply connected domains $(H_t)$ contained in the upper-half plane, onto $\mathbb{H}$. We shall denote by $K_t$ the complement of $H_t$ in $\mathbb{H}$, then for all $t > 0$, $K_t$ is a relatively compact subset of $\mathbb{H}$ and the family $(K_t)$ is increasing. For each value $\kappa > 0$, this defines a random process denoted by $\mathrm{SLE}_\kappa$ (see, e.g., [19] for more details on SLE).

In three cases, it has now been proven that $\mathrm{SLE}_\kappa$ is the scaling limit of a discrete model. Smirnov [17] proved that $\mathrm{SLE}_6$ (which is one of the processes we will focus on in the present paper) is the scaling limit of critical site percolation interfaces on the triangular grid, and Lawler, Schramm and Werner [13] have proved that $\mathrm{SLE}_2$ and $\mathrm{SLE}_8$ are the respective scaling limits of planar loop-erased random walks and uniform Peano curves. In fact, we will use Smirnov's result as a key argument in the present paper.

It is natural to study the geometry of $\mathrm{SLE}_\kappa$, and in particular, its dependence on $\kappa$. It is known (see [15] for $\kappa \neq 8$ and [13] for $\kappa = 8$) that there almost surely exists a continuous curve $\gamma : [0, \infty) \to \bar{\mathbb{H}}$ (called the *trace* of the SLE) that generates $K_t$, in the following sense: $H_t$ is the (unique) unbounded connected component of $\mathbb{H} \setminus \gamma([0, t])$. Furthermore (see [15]), $\gamma$ is a simple curve when $\kappa \leq 4$, and it is a space-filling curve when $\kappa \geq 8$.

It is possible, for each $x \in \mathbb{H}$, to evaluate the asymptotics when $\varepsilon \to 0$ of the probability that $\gamma$ intersects the disk of radius $\varepsilon$ around $x$. When $\kappa < 8$, this probability decays like $\varepsilon^\alpha$ for some $\alpha = \alpha(\kappa) > 0$. This (loosely speaking) shows that the expected number of balls of radius $\varepsilon$ needed to cover $\gamma[0, 1]$ (say) is of the order of $\varepsilon^{-2+\alpha}$, and implies that the Hausdorff dimension of $\gamma$ is not larger than $2 - \alpha$. Rohde and Schramm [15] used this strategy to show that almost surely the Hausdorff dimension of the $\mathrm{SLE}_\kappa$ trace is not larger than $1 + \kappa/8$ when $\kappa \leq 8$.

This exponent $\alpha$ and various other exponents describing exceptional subsets of $\gamma$ are closely related to critical exponents that describe the behavior near the critical point of some functionals of the related statistical physics model. Actually, in the physics literature, the derivation of the exponent is often announced in terms of (almost sure) fractal dimension, thereby omitting to prove the lower bound on the dimension. Indeed, it may *a priori* be the case that the value $\varepsilon^{-2+\alpha}$ is due to exceptional realizations of SLE$_\kappa$ with exceptionally many visited balls of radius $\varepsilon$, while "typical" realizations of SLE$_\kappa$ meet many fewer disks. One usual way to exclude such a possibility and to prove that $2 - \alpha$ corresponds to the almost sure dimension of a random fractal is to estimate second moments, that is, given *two* balls of radius $\varepsilon$, to estimate the probability that the SLE trace intersects both of them.

It is conjectured that for all $\kappa \in [0, 8]$, the Hausdorff dimension of the trace of SLE$_\kappa$ is indeed almost surely $1 + \kappa/8$. Up to the present paper, this is known to hold for $\kappa = 8/3$ for reasons that will be described below. We prove that it is the case for $\kappa = 6$:

THEOREM 1. *Almost surely, the dimension of the* SLE$_6$ *trace is* $7/4$.

Note that the discrete analog of this theorem in terms of percolation is an open problem, while it is known that the expected number of steps of a discrete exploration process is $N^{7/4}$ (cf. [18] for further reference).

Another natural object is the *boundary* of an SLE, namely $\partial K_t \cap \mathbb{H}$. For $\kappa \leq 4$, since $\gamma$ is a simple curve, the boundary of the SLE is the SLE itself; for $\kappa > 4$, it is a strict subset of the trace, and its dimension is conjectured to be $1 + 2/\kappa$. Again, the first-moment estimate is known to hold for all $\kappa$, but the only value of $\kappa > 4$ for which the dimension is known rigorously is $\kappa = 6$:

THEOREM 2 ([6]). *Almost surely, the dimension of the* SLE$_6$ *boundary is* $4/3$.

It is known that SLE$_6$ is closely related to planar Brownian motion, so that this theorem is equivalent to the same statement for the exterior boundary of a Brownian path. It was first conjectured by Mandelbrot that the fractal dimension of the boundary should be $4/3$; the first mathematical proof is due to Lawler, Schramm and Werner (cf. [6] for a review) and goes as follows.

First, note that to each point of the Brownian path, two independent Brownian motions can be associated (the past and the future), and that this point is on the boundary of the complete path iff the union of these two processes does not disconnect it from infinity. This remark provides a relation between the dimension of the boundary and the nondisconnection exponent for two paths. It is then necessary to compute the value of this exponent, and this requires a long and very technical proof. In particular, it uses the fact that the Brownian intersection exponents are analytic [10] and sharp estimates for the probabilities of nondisconnection events. (These estimates, up to the value of the exponents, were obtained earlier by Lawler in a series of clever and technical papers.)

It is conjectured (see [15] for a discussion) that the boundary of $\mathrm{SLE}_\kappa$, $\kappa > 4$, is very similar to the trace of $\mathrm{SLE}_{16/\kappa}$, and a precise statement of this *duality* is known for $\kappa = 6$ (see [11]). This and Theorem 2 provide the dimension of $\mathrm{SLE}_{8/3}$, namely: With probability 1, the dimension of the $\mathrm{SLE}_{8/3}$ trace is $4/3$.

In the present paper, we will reprove, without using the relation to planar Brownian motion, that the dimension of the outer frontier of $\mathrm{SLE}_6$ is almost surely $4/3$. Combining this with the previously mentioned universality arguments, this implies also that the dimension of the $\mathrm{SLE}_{8/3}$ trace and that of the outer frontier of planar Brownian motion are almost surely $4/3$ and gives a shorter proof of these results. We should also mention here that $\mathrm{SLE}_{8/3}$ is the natural candidate for the scaling limit of self-avoiding walks [12] and therefore also an interesting object.

Theorem 1 can be related to the dimension of *pioneer points* on a Brownian path (i.e., points $B_t$ that are on the boundary at time $t$): It is known [8] that the set of pioneer points has dimension $7/4$, the same as the $\mathrm{SLE}_6$ trace, and this is not surprising since they play similar roles. However, it can be proved that they are different. (Note, e.g., that Brownian motion can enter its past hull and the SLE trace cannot.)

The method described here cannot be extended directly to other values of $\kappa$. Indeed, two properties that are specific to $\mathrm{SLE}_6$ are used, namely the chordal/radial equivalence (in the computation of the hitting probabilities) and the locality property (in the derivation of second moments). It should be possible to obtain second moments using only the Markov property (at the cost of a more technical proof); however, the derivation of the hitting probabilities will need a different approach.

It is also possible to compute the dimension of exceptional time-sets. This is in fact easier than for subsets of the upper-half plane, since the distortion of space due to the past does not influence the probability estimates, and this makes it possible to compute dimensions for every $\kappa \geq 0$. In the last section we compute the dimension of the set of boundary times and that of the set of cut-times [i.e., times $t$ such that $\gamma(t)$ is, resp., a boundary point or a cut-point of $K$]. In particular, we prove the following:

THEOREM 3. *Let $(K_t)$ be an $\mathrm{SLE}_\kappa$ for $\kappa < 8$. Then, almost surely, $K_1$ has cut-points.*

A natural question is then the way subsets of the time interval are mapped into the upper-half plane by the trace of an SLE process. It is known [4] that Brownian motion in the plane doubles Hausdorff dimensions (i.e., if $B$ is a planar Brownian motion, then with probability 1, for all Borelian subsets $A$ of the real line, the Hausdorff dimension of $\{B_t, t \in A\}$ is equal to twice that of $A$). Such a result does not seem to hold in the general case for an SLE process.

**1. Ingredients.** We provide in this section several estimates and tools which will be needed in the subsequent proofs, but are also (maybe) of more general interest.

1.1. *Hausdorff dimension of random sets.* We will use the following result to derive the value of Hausdorff dimensions from the values of exponents. It is stated here in dimension $d \geq 1$, but we will use it only for $d = 1$ (for time sets) or $d = 2$ (for subsets of the complex plane).

Suppose that $\lambda$ denotes the Lebesgue measure in $[0,1]^d$. Let $(C_\varepsilon)_{\varepsilon > 0}$ be a family of random Borelian subsets of the cube $[0,1]^d$. Assume that for $\varepsilon < \varepsilon'$ we have $C_\varepsilon \subseteq C_{\varepsilon'}$, and let $C = \bigcap C_\varepsilon$. Define the following conditions (where $f \asymp g$ means that there exist positive numbers $c_-$ and $c_+$ such that $c_- g \leq f \leq c_+ g$, and where the constants do not depend on $\varepsilon$, $x$ or $y$):

1. For all $x \in [0,1]^d$,
$$P(x \in C_\varepsilon) \asymp \varepsilon^s.$$

2. There exists $c > 0$ such that, for all $x \in [0,1]^d$ and $\varepsilon$,
$$P(\lambda(C_\varepsilon \cap \mathcal{B}(x, \varepsilon)) > c\varepsilon^d | x \in C_\varepsilon) \geq c > 0.$$

3. There exists $c > 0$ such that, for all $x, y \in [0,1]^d$ and $\varepsilon$,
$$P(\{x, y\} \subset C_\varepsilon) \leq c\varepsilon^{2s}|x - y|^{-s}.$$

PROPOSITION 1. (i) *If conditions* 1 *and* 2 *hold, then a.s.* $\dim_H(C) \leq d - s$.

(ii) *If conditions* 1 *and* 3 *hold, then with positive probability* $\dim_H(C) \geq d - s$.

PROOF. A detailed proof of this proposition can be found in [3], Theorem 2. The outline goes as follows. First, if conditions 1 and 2 hold, they provide an upper bound on the expected number of balls of radius $\varepsilon$ needed to cover $C_\varepsilon$, hence $C$. By Borel–Cantelli, this gives an upper bound on the Minkowski dimension of $C$, which is valid with probability 1.

To derive a lower bound, one introduces the random measures $\mu_\varepsilon$ having density $\varepsilon^{-s} \mathbb{1}_{C_\varepsilon}$ with respect to the Lebesgue measure in $[0,1]^d$. If conditions 1 and 3 hold, with positive probability it is possible to extract a subsequence $\mu_{\varepsilon_k}$ converging to some measure $\mu$ supported on $C$, and to prove that with positive probability $\mu$ is a Frostman measure with dimension $d - s$, which implies that the Hausdorff dimension of the support of $\mu$ is at least $d - s$. □

Each time we will derive an almost sure Hausdorff dimension, we will in fact check these three conditions and use a zero–one law to conclude.

REMARK 1. A similar proposition can be found in [5], stated in a discrete setup in which condition 2 does not appear. Indeed, in most cases, this condition is a direct consequence of condition 1 and the definition of $C_\varepsilon$ (for instance, if $C_\varepsilon$ is a union of balls of radius $\varepsilon$).

1.2. *An estimate for diffusions.* We will need estimates for stochastic flows in an interval, that we now state and prove. For background on this topic see, for instance, [1].

Let $(X_t)$ be the diffusion process on the interval $I = [-1, 1]$ defined by the following stochastic differential equation:

$$dX_t = \sigma \, dB_t + f(X_t) \, dt, \tag{2}$$

where $\sigma > 0$ and $f$ is a given smooth function satisfying $f' \leq -a < 0$ and

$$-f(-1 + x) \sim f(1 - x) \sim -C_0 \cdot x^{-1}, \tag{3}$$

$$f'(-1 + x) \sim f'(1 - x) \sim -C_1 \cdot x^{-2}, \tag{4}$$

$$-f''(-1 + x) \sim f''(1 - x) \sim -C_2 \cdot x^{-3}, \tag{5}$$

as $x \to 0+$ for some positive constants $C_1$, $C_2$, $C_3$.

Let $(g_t)$ be the stochastic flow associated to this stochastic differential equation; that is, $(g_t)_{t \geq 0}$ is the family of random functions from $I$ to itself such that $g_t(x)$ is the value at time $t$ of the solution of (2) starting from $x$ at $t = 0$. Note that $X$ is absorbed on the points 1 and $-1$. This implies that, with probability 1, for all $t > 0$, there is an interval $I_t \subset I$ such that

$$g_t(I) = \{-1, 1\} \cup I_t.$$

We want to estimate the length $l_t$ of $I_t$. Consider the following family of partial differential equations, indexed by $b \geq 0$:

$$(E_b) \qquad \dot{h}(t, x) = \frac{\sigma^2}{2} h''(t, x) + f(x) h'(t, x) - b f'(h) h(t, x).$$

Assume that for each $b \geq 0$, $(E_b)$ has a positive solution $h_b(t, x)$ satisfying

$$h_b(t, x) \asymp [(1 + x)(1 - x)]^{q(b)} e^{-\lambda(b) t}. \tag{6}$$

It is then possible, using the Feynman–Kac formula (following exactly [8]), to prove that if $b > 0$,

$$E((g'_t(x))^b) \asymp e^{-\lambda(b)t}[(1+x)(1-x)]^{q(b)} \tag{7}$$

whenever $t \geq 1$ (where as usual we let $g'_t(x) = 0$ if the path starting from $x$ is absorbed by the boundary before time $t$). For all $x$, let

$$\tau_x = \mathrm{Inf}\{t : g_t(x) \in \{-1, 1\}\} = \mathrm{Inf}\{t : g'_t(x) = 0\}.$$

It will be sufficient for our purposes to provide an estimate in the case $\tau_0 > t$, that is, when we condition the path from 0 not to exit the interval by time $t$.

LEMMA 1. *In the previous setup,*

$$\forall b > 0, \qquad E(l_t^b \mathbb{1}_{\tau_0 > t}) \asymp e^{-\lambda(b)t}.$$

Note that this type of result does not seem to be standard in the literature on diffusions. The natural way to obtain estimates on the length of $I_t$ is to use Jensen's inequality, and, depending on the value of $b$, it can give a lower bound (if $b < 1$) or an upper bound (if $b > 1$) of the right form. Another way to obtain a lower bound is given in [8], and consists in computing the length of the image of a small interval around 0, thus giving a lower bound in terms of $g'(0)$ which is valid for all $b > 0$. Hence, all that needs to be done to complete the proof is to derive the upper bound in the case $b < 1$.

PROOF OF LEMMA 1. A rough idea of the proof is as follows: Write the length $l_t$ as the integral of $g'_t$ over $I$, and then obtain a uniform upper bound on $g'_t(x)$. Roughly speaking, two cases need to be considered (and they will be treated quite differently):

  (i) If $g_s(0)$ stays away from the boundary for $s \leq t$, then so does $g_s(x)$ for each $x$; in this case $g''_t(x)$ is bounded above and it will be possible to compare $g'_t(x)$ to $g'_t(0)$ and use (7).
  (ii) If $g_s(0)$ comes close to the boundary for $s \leq t$, then so does $g_s(x)$ for each $x$ and in this case $g'_t(x)$ becomes very small.

In the proof, $C$ will stand for a generic positive constant, the value of which may change from line to line.

We first consider the first case. The definition of $g$ implies that, for all $x \in I$,

$$g'_t(x) = \exp\left[\int_0^t f'(g_s(x))\, ds\right], \tag{8}$$

and differentiating this with respect to $x$ leads to

$$\frac{g_t''(x)}{g_t'(x)} = \int_0^t g_s'(x) f''(g_s(x)) \, ds. \tag{9}$$

Moreover, since $f'$ is bounded by $-a < 0$, (8) also proves that almost surely, for all $t > 0$ and for all $x \in I$,

$$g_t'(x) \leq e^{-at} \tag{10}$$

and in particular $l_t \leq 2e^{-at}$.

Let $\alpha > 0$ and $J_s = [-1 + \alpha e^{-as/4}, 1 - \alpha e^{-as/4}]$: If for all $s > 0$, $g_s(x) \in J_s$, then condition (5) leads to $|f''(g_s(x))| \leq C_2 \alpha^{-3} e^{3as/4}$, hence

$$\left| \frac{g_t''(x)}{g_t'(x)} \right| \leq \int_0^t C_2 \alpha^{-3} e^{-as} e^{3as/4} \, ds \leq 4C_2 a^{-1} \alpha^{-3}.$$

Assume that for all $s \in [0, t]$, $g_s(0) \in J_s$ (so that the previous estimate applies for $x = 0$). For all $x \in (-1, 0)$ such that $\tau_x > t$, write

$$\int_0^t f'(g_s(x)) \, ds = \int_0^t f'(g_s(x)) \mathbb{1}_{g_s(x) \in J_s} \, ds + \int_0^t f'(g_s(x)) \mathbb{1}_{g_s(x) \notin J_s} \, ds.$$

In the first integral, integrating $f''$ over $[g_s(x), g_s(0)]$ (which is a subset of $J_s$) and using (5) shows that

$$|f'(g_s(x)) - f'(g_s(0))| \leq Ce^{-as}(\alpha e^{-as/4})^{-3} = C\alpha^{-3} e^{-as/4}. \tag{11}$$

In the second one, since $g_s$ is monotonous, $g_s(x)$ can only be in $[-1, -1 + \alpha e^{-as/4}]$, on which $f'$ is negative and increasing. Hence, $f'(g_s(x)) \leq f'(-1 + \alpha e^{-as/4})$, and integrating $f''$ between $-1 + \alpha e^{-as/4}$ and $g_s(0)$ as previously leads to

$$|f'(-1 + \alpha e^{-as/4}) - f'(g_s(0))| \leq C\alpha^{-3} e^{-as/4}.$$

In both cases we finally obtain

$$f'(g_s(x)) \leq f'(g_s(0)) + Ce^{-as}(\alpha e^{-as/4})^{-3},$$

and integrating over $s \in [0, t]$ then proves that

$$g_t'(x) \leq \exp\left[C\alpha^{-3} + \int_0^t f'(g_s(0)) \, ds\right] \leq Kg_t'(0).$$

A similar computation shows that this also holds for $x \in (0, 1)$. Integrating this inequality leads to $l_t \leq 2Kg_t'(0)$, hence to the desired conclusion—on the event $\{\forall s \in [0, t], g_s(0) \in J_s\}$. Together with our estimate on $|g'|$ [given by (7)], this can be rewritten as

$$E[l_t^b \mathbb{1}_{\forall s < t, g_s(0) \in J_s}] \leq Ce^{-\nu t}.$$

Note that the same upper bound would apply if we replaced everywhere 0 by any given point $x$; we would even get an additional factor $(1-x^2)^q$, where $q = q(b)$ is the same as in (7). It would also apply if we considered $g_t(I_0)$ for some sub-interval $I_0$ of $[-1,1]$, in which case $g'$ would only be integrated on $I_0$ and we would get yet another multiplicative factor of the form $l(I_0)^b$.

Now let $t_1 < t_2$, and let $\mathcal{F}_{t_1}$ be the $\sigma$-field generated by the diffusion up to time $t_1$. Applying the Markov property at time $t_1$ and then the previous argument on the time interval $[t_1, t_2]$, we obtain the following estimate:

$$(12) \qquad E[l_{t_2}^b \mathbb{1}_{\forall s \in (t_1, t_2), g_s(0) \in J_s} | \mathcal{F}_{t_1}] \leq C(1 - g_{t_1}(0)^2)^q l_{t_1}^b e^{-\nu(t_2 - t_1)}.$$

We now have to consider the case where $g_s(0)$ exits $J_s$ (and this will happen in particular for small values of $s$, for which $J_s$ can even be empty if $\alpha$ is large enough). For this, we shall count the number of times it does it and use the previous estimate (12) between those times.

First, assume that for some stopping time $t > 0$ we have $g_t(0) = -1 + \alpha e^{-at/4}$ (i.e., it is on the boundary of $J_t$). Assume that $\alpha$ was chosen large enough so that the whole image of $g_t$ (except maybe the point $\{1\}$) is contained in $[-1, -1 + 2\alpha e^{-at/4}]$ [this is possible by (10)], and assume also that $t$ is sufficiently large so that $2\alpha e^{-at/4} < 1$. Let $(\varphi_s)$ denote the stochastic flow started at time $t$; for each $x$ satisfying $\tau_x > t$, we have $g_{t+s}(x) = \varphi_s(g_t(x))$ (but $(\varphi_s)$ is defined on the whole interval $[-1,1]$).

Up to the first time $\tau$ when $\varphi_\tau(-1 + 2\alpha e^{-at/4}) = 0$, $f'$ is increasing on the image of $[-1, -1 + 2\alpha e^{-at/4}]$ by $\varphi$, so we obtain $l_{t+s} \leq l_t \varphi'_s(-1 + 2\alpha e^{-at/4})$. But without a lower bound on $\tau$, we cannot apply (7) to $\varphi$ at time $\tau \wedge 1$ directly.

Instead, make the following remark. Let $u > 0$, and let $s$ be a stopping time for which $g_s(0) = -1 + u$; let $\sigma$ be the first time after $s$ for which $g_\sigma(0) = -1 + 2u$. Scaling shows that, with probability at least $\eta > 0$ (independent of $u$), we have $\sigma - s \geq u^2$. In particular, this implies that with positive probability, $g'_\sigma(0) \leq c g'_s(0)$ for some $c < 1$. Since $g'$ is decreasing anyway, we obtain $E(g'_\sigma(0)^b | \mathcal{F}_s) \leq c' g'_s(0)^b$ with $c' < 1$. We can then apply the result for $u = 2\alpha e^{-at/4}$; then $u = 2^k \alpha e^{-at/4}$ with $k \geq 2$, as long as $2^{k+1} \alpha e^{-at/4} < 1$. The number of steps we can perform is of the order of $\log_2(e^{-at/4})$, that is, linear in $t$, thus providing the estimate

$$E(l_\tau^b | \mathcal{F}_t, g_t(0) = -1 + \alpha e^{-at/4}) \leq e^{-\eta' t} l_t^b$$

with $\eta' > 0$. Combining that with (7) at time 1, we finally obtain

$$(13) \qquad E(l_{t+(\tau \wedge 1)}^b | \mathcal{F}_t, g_t(0) = -1 + \alpha e^{-at/4}) \leq C l_t^b e^{-\eta'' t/4}$$

with $\eta'' = (qa/4) \wedge \eta' > 0$.

Now, let $E_{t;k;n_1,\ldots,n_k}$ be the event that before time $t$, $g_s(0)$ exits $J_s$ on each of the $[n_i, n_{i+1}]$ and none of the others (so that roughly speaking it gets close to the boundary $k$ times). For each $i$, let $t_i$ be the first time in $[n_i, n_i+1]$ at which $g_t(0)$ is on the boundary of $J_t$. We can apply the previous estimates (12) and (13) at times $t_i$, leading to

$$E(l^b_{n_i+2}; E_{t;k;n_1,\ldots,n_k}|\mathcal{F}_{n_i}) \leq C l^b_{n_i} e^{-\eta'' n_i}.$$

We can always assume that, for each $i$, $n_{i+1} \geq n_i + 2$. Recursive application of the Markov property leads to

$$E[l^b_t \mathbb{1}_{\tau_0 > t}; E_{t;k;n_1,\ldots,n_k}] \leq C \cdot C^k e^{-\nu L_t} \prod_{i=1}^{k} e^{-\eta'' n_i},$$

where $L_t = t - \sum(n_i + 2 - n_i) = t - 2k$; so that, replacing $C$ by $Ce^{2\nu}$, we obtain

$$E[l^b_t \mathbb{1}_{\tau_0 > t}; E_{t;k;n_1,\ldots,n_k}] \leq C \cdot C^k e^{-\nu t} \prod_{i=1}^{k} e^{-\eta'' n_i}.$$

Summing over all possible values for the $n_i$ (note that we always have $n_i \geq i$) and over $k \geq 0$, we obtain

$$E[l^b_t \mathbb{1}_{\tau_0 > t}] \leq Ce^{-\nu t} \sum_{k=0}^{\infty} C^k \prod_{i=1}^{k} \frac{e^{-\eta'' i}}{1 - e^{-\eta''}} \leq Ce^{-\nu t} \sum_{k=0}^{\infty} C^k e^{-\eta'' k^2/2}.$$

This last sum being finite, we get the result we wanted. □

LEMMA 2. *In the same setup, the probability that a given point $x$ survives up to time $t > 0$ is*

$$P(\tau_x > t) = P(g'_t(x) > 0) \asymp e^{-\lambda(0)t}.$$

PROOF. We know that $E(h_0(0, g_t(x))) = h_0(t, x) \asymp e^{-\lambda(0)t}$. On the other hand, since $h_0$ is bounded, we have $E(h_0(0, g_t(x))) \leq \|h_0(0, \cdot)\|_\infty P(\tau_x > t)$. Hence

$$P(\tau_x > t) \geq \frac{ce^{-\lambda(0)t}}{\|h_0(0, \cdot)\|_\infty}.$$

Conversely, consider the distribution of $g_1(x)$. It is easy to see that, except for Dirac masses at $-1$ and $1$, it has a bounded density $p_1$ with respect to the Lebesgue measure. Since $h$ is positive, we know that $-\lambda(0)$ is the largest eigenvalue of the generator of the diffusion, and that it is simple; hence, $\|p_t\|_2 \leq \|p_1\|_2 \exp(-(t-1)\lambda(0))$. It is then a direct application of the Cauchy–Schwarz inequality to see that $\|p_t\|_1 \leq Ce^{-\lambda(0)t}$, and since we have $\|p_t\|_1 = P(\tau_x > t)$, this completes the proof of the lemma. □

**2. Dimension of the trace of SLE$_6$.**

2.1. *Construction of the trace.* Let $K$ be a chordal SLE in the upper-half plane and let $C$ be the intersection of its trace with the square $[-1,1]\times[1,3]$. In order to apply Proposition 1, introduce

$$C_\varepsilon = \{z \in [-1,1]\times[1,3] : d(z,C) \leq \varepsilon\}.$$

Since $C$ is a compact set, we have $C = \bigcap C_\varepsilon$. Moreover, we make the following remark: Let $z$ be some point in $[-1,1]\times[1,3]$, $\varepsilon > 0$, and assume that $z$ is at distance greater than $\varepsilon$ from the boundary of the square. Let $T_{\mathcal{B}(z,\varepsilon)}$ be the hitting time defined as usual as

$$T_{\mathcal{B}(z,\varepsilon)} := \mathrm{Inf}\{t : K_t \cap \mathcal{B}(z,\varepsilon) \neq \varnothing\}.$$

Then, we have the following equivalence:

(14) $\qquad z \in C_\varepsilon \quad \Longleftrightarrow \quad \mathcal{B}(z,\varepsilon) \not\subset K_{T_{\mathcal{B}(z,\varepsilon)}}.$

We call the second part of the equivalence *nondisconnection*. Indeed, the condition is equivalent to the fact that $K_{T_{\mathcal{B}(z,\varepsilon)}}$ does not disconnect $z$ from $\infty$. Note the similarity with the definition of Brownian pioneer points [5].

2.2. *The (non)disconnection exponent.* The proofs in this section rely on the equivalence between chordal and radial SLE for $\kappa = 6$ that have been proven in [8]. More precisely, there are two versions of SLE in the unit disk. The first one (chordal SLE in the disk) is obtained by mapping chordal SLE in the upper-half plane to the disk by a conformal map, so that it grows toward a point on the unit circle. The second version is called *radial* SLE, and it corresponds to the case where $K$ grows toward 0 instead of a point on the boundary. It is defined by the following PDE [if $(\tilde{g}_t)$ is the corresponding family of conformal maps]:

$$\partial_t \tilde{g}_t(z) = \tilde{g}_t(z) \frac{\tilde{g}_t(z) + \beta_t}{\tilde{g}_t(z) - \beta_t},$$

where $\beta_t = e^{i\sqrt{\kappa} t}$ is a time-scaled Brownian motion on the unit circle.

Chordal/radial equivalence is stated as follows. Let $(K_t)$ be a *chordal* SLE$_6$ in the unit disk, starting at 1 and aiming at $-1$, and let $(\tilde{K}_t)$ be a *radial* SLE$_6$ in the unit disk, starting from 1 and aiming at 0. Let $T$ (resp. $\tilde{T}$) be the first time when $K$ (resp. $\tilde{K}$) separates $-1$ from 0. Then, $K_{T^-}$ and $\tilde{K}_{\tilde{T}^-}$ have the same law, and so do $(K_{t \wedge T})_{t>0}$ and $(\tilde{K}_{t \wedge \tilde{T}})_{t>0}$ up to a (random) time change. For complete references about this, see [8]. Note that this is specific to the case $\kappa = 6$.

PROPOSITION 2. *Let $(K_t)$ be a chordal SLE$_6$ in the unit disk, starting from 1 and growing toward $-1$, and let $T_r$ be the first time when $K_t$ hits the ball with radius $r$ centered at 0. Then $K_{T_r}$ disconnects this disk from the unit circle if and only if $\mathcal{B}(0,r) \subset K_{T_r}$, and as $r$ tends to 0,*

$$p(r) := P(\mathcal{B}(0,r) \not\subset K_{T_r}) \asymp r^{1/4}.$$

PROOF. This estimate is similar to Theorem 3.1 in [8], of which it is the natural counterpart in the case $b = 0$. Let $K'$ be a *radial* $\mathrm{SLE}_6$ in the unit disk, aiming at 0, and let $T'_r$ be the first time when it reaches $\mathcal{B}(0,r)$. The chordal/radial equivalence shows that $p(r)$ is equal to the probability that $K'_{T'_r}$ does not disconnect $\mathcal{B}(0,r)$ from $-1$, that is, the probability that $-1 \notin K'_{T'_r}$.

Let $W_t = e^{i\sqrt{6}B_t}$ be the (time-scaled) Brownian motion on $\partial \mathbb{U}$ driving $(K'_t)$ [where $(B_t)$ is a standard Brownian motion on $\mathbb{R}$], and let $Y_t$ be the continuous determination of the argument of $g_t(-1)/W_t$ starting at $\pi$. $Y_t$ is well defined as long as $K'$ does not reach $-1$. Loewner's differential equation and Itô's formula show that

$$dY_t = \sqrt{6}\, dB_t + \cotg(Y_t/2)\, dt,$$

that is, $(Y_t)$ is a diffusion process with diffusion $\sqrt{6}$ and drift $\cotg(\cdot/2)$, absorbed by $\{0, 2\pi\}$ when $-1$ is absorbed by $K'_t$. Straightforward calculations prove that $f_t(x) = e^{-t/4}(\sin y/2)^{1/3}$ satisfies $\partial_t f_t = Lf_t = -\frac{1}{4}f_t$; we can now apply Lemma 2 and obtain

(15) $$P(-1 \notin K'_t) \asymp e^{-t/4}.$$

But Köbe's distortion theorem [14] states that, if $r(t) = d(0, K_t)$, then

$$\frac{e^{-t}}{4} \leq r(t) \leq e^{-t},$$

which, combined with estimate (15), proves the proposition. (Details of the last step are the same as in [8].) □

COROLLARY 1. *Fix $\eta > 0$, and let $\mathcal{B} = \mathcal{B}(z, r)$ be some disk contained in $\mathbb{U}$, where $|z| < 1 - 2\eta$ and $r < \eta$; let $(K_t)$ be a chordal $\mathrm{SLE}_6$ in the unit disk, starting from 1 and aiming at $-1$. If $T_\mathcal{B}$ denotes the first time when $K_t$ reaches $\mathcal{B}$, then the probability $p(\mathcal{B})$ that $K_{T_\mathcal{B}}$ does not disconnect $\mathcal{B}$ from $-1$ satisfies*

$$p(\mathcal{B}) \asymp r^{1/4},$$

*where the implicit constants depend only on $\eta$.*

PROOF. There exists exactly one Möbius transform $\Phi : \mathbb{U} \to \mathbb{U}$ mapping 1 to itself and $\mathcal{B}$ to a disk centered at 0. The radius of $\Phi(\mathcal{B})$ is then

$$\rho(z, r) = \frac{(1 + r^2 - |z|^2) - \sqrt{(1 + r^2 - |z|^2)^2 - 4r^2}}{2r} \asymp r.$$

$\Phi(K)$ is then a chordal SLE in the disk starting from 1 and aiming at $\Phi(-1)$. Moreover, $|\Phi(-1) - 1|$ is bounded away from 0 by a constant. The proof of

Proposition 2 can then be adapted (only changing the position of the end point) to show that

$$p(\mathcal{B}) \asymp p(\rho(z,r)) \asymp p(r) \asymp r^{-1/4},$$

with constants depending only on $\eta$. □

It is then easy, by mapping the disk to the upper-half plane and using (14), to turn this corollary into the first condition of Proposition 1, that is,

(16) $\qquad \forall z \in [-1,1] \times [1,3], \qquad P(z \in C_\varepsilon) \asymp \varepsilon^{1/4}.$

It then follows from the definition of $C_\varepsilon$ that condition 2 holds: If $z \in C_\varepsilon$, let $z' \in C$ such that $|z - z'| = \varepsilon$ (which exists by a compactness argument); then the disk with diameter $[zz']$ is contained in $\mathcal{B}(z,\varepsilon) \cap C_\varepsilon$ and it has area $\pi \varepsilon^2 / 4$.

2.3. *Percolation and second moments.* We now turn our attention to the derivation of second moments for the hitting probability of disks by the SLE$_6$ trace, namely condition 3 in Proposition 1. Again we will make strong use of the fact that we are in the case $\kappa = 6$, and in fact the decay of correlations we obtain is a consequence of the locality property of SLE$_6$. It has been proved [17, 18] that the exploration process of critical percolation on the triangular lattice converges to the SLE$_6$ trace; in particular, consider critical percolation on a discretization of the upper-half plane with mesh $\delta > 0$ and the usual boundary conditions (i.e., wired on $[0, +\infty)$ and free on $(-\infty, 0)$). Then the probability that the discrete exploration process $\gamma_\delta$ hits the ball $\mathcal{B}(i, \varepsilon)$ satisfies

(17) $\qquad P(\gamma_\delta \cap \mathcal{B}(i,\varepsilon) \neq \varnothing) \underset{\delta \to 0}{\longrightarrow} P(i \in C_\varepsilon) \asymp \varepsilon^{1/4}.$

But the fact that the discrete exploration process touches this disk is equivalent to the existence of both a closed path connecting the disk to $[0, +\infty)$ and an open path connecting the disk to $(-\infty, 0)$. Applying the results in [18], this leads to the following:

COROLLARY 2. *Let $A_\varepsilon$ be the annulus centered at 0, with radii $\varepsilon$ and 1. For all $\delta > 0$, consider critical site-percolation on the intersection of $A_\varepsilon$ with the triangular lattice of mesh $\delta$. Let $p(\varepsilon, \delta)$ be the probability that $\mathcal{C}(0, \varepsilon)$ is connected to $\mathcal{C}(0, 1)$ by both a path of open vertices and a path of closed vertices in $A_\varepsilon$. Then, as $\delta$ tends to 0, $p(\varepsilon, \delta)$ converges to some $p(\varepsilon)$ satisfying*

$$p(\varepsilon) \asymp \varepsilon^{1/4}.$$

Note that this says nothing about the speed of convergence, and hence does not provide useful estimates for the probability of the discrete event itself; but it is sufficient for our purpose here.

Now fix $z, z' \in [-1, 1] \times [1, 3]$ and $\varepsilon < |z - z'|/2$. Again, the probability that the $SLE_6$ trace touches both $\mathcal{B}(z, \varepsilon)$ and $\mathcal{B}(z', \varepsilon)$ can be written as the limit, as $\delta$ goes to 0, of the corresponding probability for critical site-percolation on the triangular lattice with mesh $\delta$. But this implies the following:

   (i) There exist a path of open vertices and a path of closed vertices, both connecting $\mathcal{C}(z, \varepsilon)$ to $\mathcal{C}(z, |z - z'|/2)$ inside $\mathcal{B}(z, |z - z'|/2)$.
   (ii) There exist a path of open vertices and a path of closed vertices, both connecting $\mathcal{C}(z', \varepsilon)$ to $\mathcal{C}(z', |z - z'|/2)$ inside $\mathcal{B}(z', |z - z'|/2)$.
   (iii) There exist a path of open vertices and a path of closed vertices, both connecting $\mathcal{C}((z+z')/2, |z - z'|)$ to the real axis outside $\mathcal{B}((z+z')/2, |z - z'|)$.

Those three events are independent, since they describe the behavior of pairwise disjoint sets of vertices; besides, the probability of each of them can be estimated using Corollary 2 and converges, as $\delta \to 0$ and up to universal multiplicative constants, respectively, to $(\varepsilon/d)^{1/4}$, $(\varepsilon/d)^{1/4}$ and $d^{1/4}$, with $d = |z - z'|$. Hence, letting $\delta$ go to 0, we obtain the following estimate:

$$(18) \quad P(\{z, z'\} \subset C - \varepsilon) \leq C \left( \frac{\varepsilon}{|z - z'|} \right)^{1/2} |z - z'|^{1/4} = C \frac{\varepsilon^{1/2}}{|z - z'|^{1/4}},$$

which is exactly condition 3 in Proposition 1 with $s = 1/4$, as we wanted.

2.4. *Conclusion.* It is now possible to apply Proposition 1 with $s = 1/4$: We obtain

$$P(\dim_H(C) \leq \tfrac{7}{4}) = 1, \qquad P(\dim_H(C) = \tfrac{7}{4}) > 0.$$

Now let $\mathcal{H}_\infty$ be the complete trace of $K$. Since $C \subset \mathcal{H}_\infty$, we obtain the same results for $\mathcal{H}_\infty$. Theorem 1 is then a consequence of the following:

LEMMA 3 (0–1 law for the trace). *For all $d \in [0, 2]$, we have*

$$P(\dim_H(\mathcal{H}_\infty) = d) \in \{0, 1\}.$$

PROOF. For all $n \in \mathbb{Z}$, let $D_n = \dim_H(\mathcal{H}_{2^n})$. For all $n$, we then have $D_{n+1} \geq D_n$ [because $(\mathcal{H}_t)$ is increasing] and besides $D_n$ and $D_{n+1}$ have the same law (by the scaling property). Hence, almost surely, for all $m$, $n$, we have $D_n = D_m$. Taking this to the limit gives $P(\dim_H(\mathcal{H}_\infty) = D_n) = 1$; hence the random variable $\dim_H(\mathcal{H}_\infty)$ is $\mathcal{F}_{2^n}$-measurable for all $n$. Hence it is $\mathcal{F}_{0+}$-measurable, and we know by Blumenthal's zero–one law that this $\sigma$-field is trivial. $\square$

**3. Dimension of the boundary of SLE$_6$.** In this section we adapt the previous proof to compute the Hausdorff dimension of the boundary of $K$ at some fixed time.

3.1. *The escape probability.*

PROPOSITION 3. *Let $(K_t)$ be an SLE$_6$ in the half-plane, and let $\tau_R$ be the first time it reaches radius $R$. Then, as $R$ goes to infinity,*

$$P(i \notin K_{\tau_R}) \asymp R^{-1/3}.$$

Note that the corresponding result for $P(1 \notin K_{\tau_R})$ has been derived in [8].

PROOF OF PROPOSITION 3. Fix $l > 2$ (its value will be determined later). We shall suppose that $K$ starts at $l$ instead of $0$; it is easy to see that this only changes the estimates up to a fixed constant. The idea of the proof is as follows: We will prove that the conditional probability, knowing that $i \notin K_{\tau_R}$, that $K_{\tau_R}$ does not intersect the unit disk is bounded below by a positive constant. The probability that this happens satisfies $P(\bar{\mathbb{U}} \cap K_{\tau_R} = \varnothing) \asymp R^{-1/3}$ (cf. [8]), which leads to the conclusion.

We make the following preliminary remark. For each $t \geq 0$ such that $i \notin K_t$, the intersection of $\mathcal{C}(0,2)$ with $H_t = \mathbb{H} \setminus K_t$ is a union of at most countably many arcs of positive length. Because $H_t$ is simply connected, some of these arcs separate $i$ from infinity, and the first separating arc on a continuous path from $i$ to infinity does not depend on the path. We will denote this arc by $\lambda_t^{(2)}$ and call it the *relevant arc* at radius 2 and time $t$. Define $\lambda_t^{(l)}$ similarly. It is easy to see that $(\lambda_t^{(2)})$ is nonincreasing, in the sense that, for all $0 < t < t'$, we have $\lambda_{t'}^{(2)} \subseteq \lambda_t^{(2)}$.

Introduce the following stopping times (where $S_0 = T_0 = 0$), defined inductively for $n > 0$ (where as usual we let $\operatorname{Inf} \varnothing = \infty$):

$$S_n = \operatorname{Inf}\{t > T_{n-1} : \lambda_t^{(2)} \subsetneq \lambda_{T_{n-1}}^{(2)}\},$$

$$T_n = \operatorname{Inf}\{t > S_n : \lambda_t^{(l)} \subsetneq \lambda_{S_n}^{(l)}\}.$$

(Loosely speaking, $S_n$ is the first time after $T_{n-1}$ when the process touches the circle of radius 2 and $T_n$ is the first time after $S_n$ when it returns to the circle of radius $l$, except that we are only considering the relevant parts of the circles.)

Moreover, let $T = \tau_R = \operatorname{Inf}\{t : K_t \cap \mathcal{C}(0,R) \neq \varnothing\}$. Then, almost surely, the integer $N = \operatorname{Sup}\{n : T_n < T\}$ is finite and we have

$$0 = T_0 < S_1 < T_1 < \cdots < S_N < T_N < T < \infty$$

(i.e., $K$ crosses the annulus between radii 2 and $l$ only finitely many times before reaching radius $R$). In the Brownian case, $N$ would be geometric with parameter $\log(l/2)/\log(R/2)$.

Let $E_R$ and $E'_R$ be defined as

$$E_R := \{i \notin K_T\}, \qquad E'_R := \{\bar{\mathbb{U}} \cap K_T = \varnothing\}.$$

We have to estimate $P(E_R)$; from Theorem 3.1 in [7], $P(E'_R) \asymp R^{-1/3}$, and we have $P(E_R) \geq P(E'_R)$. We shall decompose $E_R$ according to the value of $N$: We can write $P(E_R) = \sum_{n=0}^{\infty} P(E_R, N = n)$. Note that if both $E_R$ and $\{N = 0\}$ hold, then $E'_R$ holds as well.

For fixed $n$, make the following remark: If there is not disconnection before $T$, then there is not disconnection for $t$ inside any $[S_k, T_k]$, for all $k \leq n$. Apply the strong Markov property at time $S_k$, and condition $K_{S_k}$ not to contain $i$. The conditional probability of disconnection between $[S_k, T_k]$ is then at least equal to the probability that an $\mathrm{SLE}_6$ in $H_{S_k}$, starting from $\gamma(S_k)$, separates $i$ from infinity before reaching the unit circle or (the relevant part of) the circle $\mathcal{C}(0, l)$. We will prove that this last probability is bounded below by a constant which depends only on $l$.

So, let $\Omega$ be the bounded connected component of $H_{S_k} \setminus \lambda_{S_k}^{(l)}$ (i.e., the one which contains $i$), and let $\Omega'$ be the connected component of $\Omega \setminus \mathcal{C}(0,1)$ having $\gamma(S_k)$ on its boundary. Let $\partial = \partial K_{S_k} \cup \mathbb{R} \cup \mathcal{C}(0,l) \cup \mathcal{C}(0,1)$. We can write the boundary of $\Omega'$, starting from $\gamma(S_k)$ and going counterclockwise, as the union of five subsets:

(i) $\partial_1 \subset \partial$;
(ii) some arc $\partial_2$ of $\mathcal{C}(0,1) \cap H_{S_k}$ which either contains $i$ or separates it from infinity in $H_{S_k}$;
(iii) $\partial_3 \subset \partial$;
(iv) $\partial_4 = \lambda_{S_k}^{(l)}$;
(v) $\partial_5 \subset \partial$

(possibly in the opposite order, which changes nothing in what follows). If $\gamma$ touches $\partial_3$ before $\partial_2 \cup \partial_4$, then it disconnects $i$ from infinity before time $T_k$.

First, estimate the probability that $\gamma$ hits $\partial_3 \cup \partial_4$ before $\partial_2$. By the locality property, this is the probability that an $\mathrm{SLE}_6$ in $\Omega'$, going from $\gamma(S_k)$ to the common endpoint of $\partial_2$ and $\partial_3$, does the same. In turn, this is a positive non-decreasing function of the extremal distance $L_1$ between $\partial_2$ and $\partial_5$ in $\Omega'$ (cf. [7]). By construction, this extremal distance is at least equal to that between $\partial_2$ and $\lambda_{S_k}^{(2)}$, which in turn is bounded below by the extremal distance between the unit circle and $\mathcal{C}(0,2)$ in $\mathbb{H}$. Hence, the probability that $\gamma$ touches $\partial_3 \cup \partial_4$ before $\partial_2$ is bounded below by some absolute constant $\varepsilon > 0$.

Now, conditionally on the previous event, if $\gamma$ hits $\partial_3$ before $\partial_4$, then it disconnects $i$ from $\infty$ before time $T_k$. The conditional probability that this does not happen is the conditional probability to hit $\partial_4$ before $\partial_2 \cup \partial_3$, and this in turn is not greater than $\varepsilon^{-1}$ times the probability to hit $\partial_4$ before $\partial_2 \cup \partial_3$ (still conditionally to the past of the process up to time $S_k$).

Again this probability is a positive, decreasing function of the extremal distance $L_2$ between $\partial_1$ and $\partial_4$ in $\Omega'$. $L_2$ is at least equal to the extremal distance between $\lambda_{S_k}^{(2)}$ and $\partial_4$ in $\Omega'$, which is bounded below by the extremal distance between $\mathcal{C}(0,2)$ and $\mathcal{C}(0,l)$ in $\mathbb{H}$. This can be made arbitrarily large by taking $l$ large enough; so for large $l$, the probability we are considering is smaller than $\varepsilon/2$. With such a value of $l$, we thus obtain $P(i \in K_{T_k}|\mathcal{F}_{S_k}) \geq 1/2$.

So, for all $k$, the (conditional, knowing the process up to time $S_k$ and knowing $i \notin K_{S_k}$) probability that there is no disconnection between times $S_k$ and $T_k$ is bounded by $1/2$. After time $T_N$, the probability to not swallow $i$ before reaching $\mathcal{C}(0,R)$, and to do so without touching (the relevant part of) $\mathcal{C}(0,2)$ is bounded above by $cR^{-1/3}$ (applying the same reasoning between radii $l$ and $R$ instead of $2$ and $l$, resp.). Hence:

$$P(E_R) = \sum_{k=0}^{\infty} P(E_R, N=k) \leq c \sum_{k=0}^{\infty} 2^{-k} R^{-1/3} \leq 2cR^{-1/3}.$$

But on the other hand, by the results in [7], we have

$$P(E_R) \geq P(E'_R) \asymp R^{-1/3},$$

so that finally we obtain what we announced:

$$P(E_R) \asymp P(E_R, N=0) \asymp P(E'_R) \asymp R^{-1/3}. \qquad \square$$

3.2. *Exponent for $b = 1/3$.*

PROPOSITION 4. *Let $(K_t)$ be a chordal* SLE$_6$ *in the unit disk, starting from $1$ and growing toward $-1$, and let $T_r$ be the first time when $K_t$ hits the ball with radius $r$ centered at $0$. Let $L_{T_r}$ be $\pi$ times the extremal distance in $\mathbb{U} \setminus K_{T_r}$ between $\mathcal{C}(0,r)$ and $\partial \mathbb{U}$. Then, as $r$ tends to $0$,*

$$E(e^{-L_{T_r}/3}) \asymp r^{2/3}.$$

PROOF. As previously, let $K'$ be a radial SLE$_6$ in $\mathbb{U}$, starting from $1$ and aimed at $0$. Then, since all the involved events satisfy nondisconnection between $\mathcal{C}(0,r)$ and $-1$ ($L_{T_r} = \infty$ iff there is disconnection), we have

$$(19) \quad q(r) := E(e^{-L_{T_r}/3}) = E(e^{-L_{T_r}/3} \mathbb{1}_{L_{T_r} < \infty}) = E(e^{-L'_{T'_r}/3} \mathbb{1}_{-1 \notin K'_{T'_r}}).$$

We shall estimate the third term, again following the steps of the proof of Theorem 3.1 in [8]. From now on, fix $b = 1/3$ and $\nu = \nu(\kappa, b) = 2/3$; since $b < 1$, we need a separate proof here. Let $l_t$ be the Euclidean length of the arc $g_t(\partial \mathbb{U} \setminus K'_t)$. The only place in [8] where $b \geq 1$ was needed was in the derivation of

$$E(l_t^b) \asymp \exp(-\nu t). \tag{20}$$

But this is exactly what Lemma 1 shows, after suitable rescaling. $\square$

3.3. *Construction of the boundary.* Again we describe the studied set as the decreasing intersection of a family $B_\varepsilon$ of subsets of the plane. Here, let

$$B_\varepsilon = \{z \notin K_1 : d(z, K_1) < \varepsilon\}.$$

In order for $z$ to be in $B_\varepsilon$, the following must happen. First, there is some point in $\mathcal{H}$ at a distance less than $\varepsilon$ from $z$; letting $T(z, \varepsilon) = \operatorname{Inf}\{t : d(z, K_t) < \varepsilon\}$, and introducing the extremal distance $L(z, \varepsilon)$ between $\mathcal{B}(z, \varepsilon)$ and $\partial \mathbb{U}$ in $\mathbb{U} \setminus K_{T(z,\varepsilon)}$, this condition is equivalent to

$$L(z, \varepsilon) < \infty.$$

Then, the SLE after $T(z, \varepsilon)$ and up to time 1 must not disconnect $z$ from "infinity" (i.e., from $\partial \mathbb{U}$), and conditionally to $K_{T(z,\varepsilon)}$ this happens with probability of order

$$e^{-L(z,\varepsilon)/3}.$$

Proposition 4 then states that $P(z \in B_\varepsilon) \asymp \varepsilon^{2/3}$. Second moments can be obtained in the same fashion as for the trace; in this case, the relevant estimate (describing in which conditions a disk intersects the boundary of the discrete exploration process) is the following. First, two crossings of different colors must ensure that the exploration process touches the disk; then a third path, disjoint of the first two, will prevent it from closing a loop around it. Hence the following consequence of Proposition 4:

COROLLARY 3. *Let $A_\varepsilon$ be the annulus centered at 0, with radii $\varepsilon$ and 1. For all $\delta > 0$, consider critical site-percolation on the intersection of $A_\varepsilon$ with the triangular lattice of mesh $\delta$. Let $\tilde{p}(\varepsilon, \delta)$ be the probability that $\mathcal{C}(0, \varepsilon)$ is connected to $\mathcal{C}(0, 1)$ both by a path of open vertices and by two disjoint paths of closed vertices in $A_\varepsilon$. Then, as $\delta$ tends to 0, $\tilde{p}(\varepsilon, \delta)$ converges to some $\tilde{p}(\varepsilon)$ satisfying*

$$\tilde{p}(\varepsilon) \asymp \varepsilon^{2/3}.$$

The rest of the construction is the same, and we obtain sufficient estimates to apply Proposition 1, this time with $s = 2/3$. We obtain

$$P(\dim_H(\partial K_1) \leq \tfrac{4}{3}) = 1, \qquad P(\dim_H(\partial K_1) = \tfrac{4}{3}) > 0,$$

and once more we need a zero–one law.

LEMMA 4 (0–1 law for the boundary). *For all $d \in [0,2]$, we have*

$$P(\dim_H(\partial K_1) = d) \in \{0,1\}.$$

PROOF. Let $D_t = \dim_H(\partial K_t)$. As previously in the case of the trace, scaling shows that the law of $D_t$ does not depend on $t > 0$. However, here $(\partial K_t)$ is not increasing anymore, so we need another argument. Let $t, t' > 0$ and consider the boundary of $K_{t+t'}$. It has two parts, namely the "new" part $\partial_1 = \partial K_{t+t'} \setminus K_t$, and the "old" part $\partial_2 = \partial K_{t+t'} \cap K_t \subset \partial K_t$. It is clear that

$$D_{t+t'} = \dim_H(\partial_1) \vee \dim_H(\partial_2),$$

hence in particular $\dim_H(\partial_1) \leq D_{t+t'}$. Besides, conformal mapping shows that $\dim_H(\partial_1)$ has the same law as $D_{t'}$, hence the same law as $D_{t+t'}$. Hence, with probability 1, $D_{t+t'} = \dim_H(\partial_2)$.

Moreover, conformal mapping also shows that $\dim_H(\partial_2)$ is independent of $\mathcal{F}_t$. This proves that, for all $t, t' > 0$, the dimension of $\partial K_{t+t'}$ is independent of $\mathcal{F}_t$. It is then a direct application of Blumenthal's zero–one law that $D_t$ has an almost sure value. □

This concludes the proof of Theorem 1.

3.4. *Dimension of* SLE$_{8/3}$. It should be theoretically possible to apply the previous construction to other values of $\kappa$, but some of the main tools that we used (namely, the radial/chordal equivalence and the restriction property) do hold only for $\kappa = 6$, so that additional arguments would be required.

For the special value $\kappa = 8/3$, the result on the frontier of SLE$_6$ makes it possible to show that the dimension of SLE$_{8/3}$ is almost surely $4/3$. More precisely, Lawler, Schramm and Werner [11] have shown that the outer boundary of the union of eight SLE$_{8/3}$'s has the same law as that of the union of five Brownian excursions. The zero–one laws previously proved for both the trace and the boundary of SLE extends to this object: Its boundary has a.s. the same dimension as the boundary of SLE$_6$ and also a.s. the same dimension as SLE$_{8/3}$. Hence these dimensions are equal, and the result follows.

**4. Time-sets for SLE$_\kappa$.** We now turn our attention to the dimension of sets of exceptional times. Note that time corresponds to the Loewner parameterization of the trace, which is in a way not the most canonical; it is not clear, for instance, whether it behaves nicely under time-reversal. More precisely, how smoothly does the Riemann map from $\mathbb{H} \setminus \gamma([t,\infty))$ to $\mathbb{H}$ evolve as $t$ increases?

A natural question that also arises is the following. Let $A$ be some (random) subset of $[0,\infty]$, and let $\gamma(A)$ be its image by the trace of a chordal SLE in the upper-half plane. Is it possible, knowing the Hausdorff dimension of $A$, to obtain that of $\gamma(A)$? Such a relation holds for Brownian motion [4]; namely, the dimension of the image is a.s. equal to twice the dimension of $A$. It is expected that such a relation cannot hold for SLE without additional requirements on $A$; however, a few cases can be treated entirely (in the sense that both the time and space dimensions can be computed in independent ways), at least for $\kappa = 6$: the trace itself, cut-points, and the boundary.

4.1. *Boundary times.* In the previous section, we derived the dimension of the boundary of SLE$_6$. The dimension of the corresponding time-set can also be computed (and it should be noted that the following is true even for $\kappa \neq 6$).

THEOREM 4. *Let $K$ be an SLE in the upper-half plane, with $\kappa > 4$, and let $D$ be the set of boundary times in $[0,1]$, that is, the set of times $t$ such that $\gamma(t) \in \partial K_1$. Then, with probability 1,*

$$\dim_H(D) = \frac{4+\kappa}{2\kappa}.$$

PROOF. It is clearly sufficient to compute the dimension of left-boundary times, namely times $t$ such that $g_1(\gamma(t)) \in (-\infty, \beta_t)$, where $\beta$ is the process driving $K$. Introduce the sets of approximate left-boundary times between $\varepsilon$ and $a$, defined by

$$D_{\varepsilon,a} = \{t : \mathrm{Inf}(\mathbb{R} \cap g_t(K_{t+\varepsilon})) = \mathrm{Inf}(\mathbb{R} \cap g_t(K_{t+a}))\}$$

(i.e., $\gamma$ may touch the real line on the right side of $K$ but not on the left side). Let $D^a$ be the intersection of the $D_{\varepsilon,a}$ when $\varepsilon \to 0$. Scaling and the Markov property show that $P(t \in D_{\varepsilon,a})$ depends only on $\varepsilon/a$. Hence, to obtain condition 1 in Proposition 1, with $s = (\kappa-4)/2\kappa$, it is sufficient to obtain the following estimate:

LEMMA 5. *Let $(K_t)$ be a chordal SLE$_\kappa$ ($\kappa > 4$) in the upper-half plane. Then as $t$ goes to infinity,*

$$p_t := P(\mathrm{Inf}(\mathbb{R} \cap K_t) = \mathrm{Inf}(\mathbb{R} \cap K_1)) \asymp t^{(4-\kappa)/2\kappa}.$$

PROOF. First, apply the Markov property of SLE at time 1 and map the picture to the upper-half plane by $\Phi = g_1 - \beta_1$. Let $Y_0 \leq 0$ be the image of $\mathrm{Inf}(\mathbb{R} \cap K_1)$ by $\Phi$. The process $(\tilde{K}_u) = (\Phi(K_{1+u}))$ is an SLE$_\kappa$, and the probability we are interested in is then given by

$$p_t = P(Y_0 \notin \tilde{K}_{t-1}).$$

Let $(\tilde{\beta}_u)$ and $\tilde{g}_u \colon \mathbb{H} \setminus \tilde{K}_u \to \mathbb{H}$ be, respectively, the process driving $\tilde{K}$ and the associated conformal maps; let $Y_u = \tilde{g}_u(Y_0) - \tilde{\beta}_u$. It is easy to see, using Itô's formula and the definition of chordal SLE, that $Y$ satisfies the following SDE (where $B$ is a standard Brownian motion):

$$(21) \qquad dY_u = \sqrt{\kappa}\, dB_u + \frac{2}{Y_t}\, du;$$

that is, up to a linear time change, $Y$ is a Bessel process of dimension $b = 1 + 4/\kappa < 2$ starting from $Y_0$. Hence, it is known that the probability that it does not hit 0 up to time $u$ behaves like $(u/Y_0^2)^{-\nu}$, where $\nu = (\kappa - 4)/2\kappa > 0$ is the index of the process. Hence,

$$p_t \asymp t^{-\nu} E(Y_0^{2\nu}) \asymp t^{-\nu},$$

as we wanted. □

This provides the right estimate,

$$P(t \in D_{\varepsilon,a}) \asymp \left[\frac{\varepsilon}{a}\right]^s,$$

where the implicit constants depend only on $\kappa$. Notice that if $t + \varepsilon$ is in $D_{\varepsilon,a}$, then $t \in D_{2\varepsilon,a}$ (because $K_t \subset K_{t+\varepsilon}$) and even $[t, t+\varepsilon] \subset D_{2\varepsilon,a}$. This and the previous estimate provide

$$P([t, t+\varepsilon] \subset D_{2\varepsilon,a} | t \in D_{2\varepsilon,a}) \geq \frac{P(t + \varepsilon \in D_{\varepsilon,a})}{P(t \in D_{2\varepsilon,a})} \geq c > 0,$$

which is condition 2. It remains to obtain second moments, and these are given by the Markov property, as follows.

Let $x < y$ be two times in $[0,1]$. If $x$ and $y$ are in $D_{\varepsilon,a}$ with $a > y - x$, then in particular $x \in D_{\varepsilon,y-x}$ and $y \in D_{\varepsilon,a}$. By the Markov property of SLE, applied at time $y$, those two events are independent. Hence we obtain

$$P(x, y \in D_{\varepsilon,a}) \leq P(x \in D_{\varepsilon,y-x}) P(y \in D_{\varepsilon,a})$$

$$\leq C\left[\frac{\varepsilon}{y-x}\right]^s \left[\frac{\varepsilon}{a}\right]^s \leq C \frac{\varepsilon^{2s}}{(y-x)^s},$$

still with $s = (\kappa - 4)/2\kappa$. This is exactly condition 3. If $a \leq y - x$, then the events $x \in D_{\varepsilon,a}$ and $y \in D_{\varepsilon,a}$ are themselves independent and the same method applies. Hence, everything is ready to apply Proposition 1: For all $a > 0$, with positive probability,

$$\dim_H(D^a) = 1 - \frac{\kappa - 4}{2\kappa} = \frac{4 + \kappa}{2\kappa}.$$

Noticing then that $D^1 \subset D \subset D^2$ hence provides

$$P\left(\dim_H(D) = \frac{4 + \kappa}{2\kappa}\right) > 0.$$

It is then easy to apply the same proof as that of Lemma 4 and obtain a zero–one law for $\dim_H(D)$, thus completing the proof. $\square$

REMARK 2. In particular, the dimension of boundary-times is never less than $1/2$, even when $\kappa \to \infty$. Note that in this case, the dimension of the Bessel process appearing in the proof tends to 1, so the exponent $1/2$ is the same as in the usual gambler's ruin estimate.

This is not surprising since, when $\kappa$ tends to $\infty$, the trace of an $\mathrm{SLE}_\kappa$ converges, after suitable rescaling, to

$$\gamma_\infty : t \mapsto (B_t, L_t^{B_t}),$$

where $B$ is standard Brownian motion and $(L_t^x)$ denotes its local time at point $x$ (cf. [2]). In the limit, the boundary times correspond to last-passage times, which have dimension $1/2$ by a reflection argument.

4.2. *Cut-times and the existence of cut-points.* We saw in the previous sections how the dimension of the trace of SLE was related to nondisconnection exponents. Here, we follow the analogy with Brownian motion to describe cut-points on the SLE trace. Let $K$ be a chordal $\mathrm{SLE}_\kappa$ and let $C$ be the set of cut-points of $K_2$ in $K_1$ (i.e., the set of points $z \in K_1$ such that $K_2 \setminus \{z\}$ is not connected). Such a point is on the boundary of $K_1$; hence if $\gamma$ is the trace of $K$, every cut-point is on $\gamma([0,1])$. We say that $t$ is a cut-time if $\gamma(t)$ is a cut-point, and note $\tilde{C}$ the set of cut-times.

THEOREM 5.

(i) *If $0 \leq \kappa \leq 4$, then $\tilde{C} = [0,1]$ and $C = K_1$.*

(ii) *If $4 < \kappa < 8$, then with positive probability $\tilde{C}$ has Hausdorff dimension $(8 - \kappa)/4$; in particular, it is nonempty, hence $C \neq \varnothing$.*

(iii) *If $\kappa > 8$, then a.s. $\tilde{C} = \varnothing$ and $K_1$ has no cut-point.*

PROOF. (i) is a direct consequence of the fact that $\gamma$ is a simple path [15], and (iii) is proved exactly like (ii) with the usual convention that a set of negative dimension is empty. Hence, we may assume that $4 < \kappa < 8$. Again we are going to apply Proposition 1, and the proof will be very similar to that of Theorem 4.

Introduce the set of approximate cut-times between $\varepsilon$ and $a$ defined as

$$C_{\varepsilon,a} = \{t \in [0,1] : \gamma([t+\varepsilon, t+a]) \cap (K_t \cup \mathbb{R}) = \varnothing\}.$$

Define $C^a$ as the (indeed nonincreasing) intersection of the $C_{\varepsilon,a}$. By the Markov property at time $t$, it follows that $P(t \in C_{\varepsilon,a})$ does not depend on $t$. Moreover, scaling shows that it is a function of $\varepsilon/a$. Hence, to obtain condition 1 in Proposition 1 with $s = (\kappa - 4)/4$, it suffices to prove the following:

LEMMA 6. *Let $K$ be an* SLE$_\kappa$ *in the upper-half plane, starting at $x \in (0,1)$, with $\kappa > 4$. Then, when $t \to \infty$,*

$$P(\{0,1\} \cap K_t = \varnothing) \asymp t^{(4-\kappa)/4}.$$

PROOF. The proof of this lemma is very similar to that of Theorem 3.1 in [9]. Two things have to be done: first, extend this theorem to the (easier) case where $w_1 = w_2 = 0$; second, translate it back to an estimate for SLE at a fixed time. Introduce the following processes: $X_t = g_t(1) - \beta_t$, $Y_t = g_t(0) - \beta_t$, where $(\beta_t)$ is the time-scaled Brownian motion driving $K$. As was seen previously, $X$ and $Y$ satisfy the following SDEs:

$$dX_t = \frac{2}{X_t} dt + \sqrt{\kappa}\, dB_t, \qquad dY_t = \frac{2}{Y_t} dt + \sqrt{\kappa}\, dB_t,$$

where $B$ is a standard real Brownian motion. Let $L_t = X_t - Y_t$ be the length of the image interval, and let $R_t = X_t/L_t$. Tedious application of Itô's formula leads to

$$dL_t = \frac{2\, dt}{L_t R_t(1-R_t)}, \qquad dR_t = \frac{2(1-2R_t)}{L_t^2 R_t(1-R_t)} dt + \frac{\sqrt{\kappa}}{L_t} dB_t.$$

Introduce the following random time-change:

$$dt(s) = \frac{L_{t(s)}^2 R_{t(s)}(1-R_{t(s)})}{2} ds;$$

then the previous system reduces to $dL_{t(s)} = L_{t(s)}\, ds$, that is, almost surely $L_{t(s)} = e^s$ and, letting $Z_s = R_{t(s)}$,

(22) $$dZ_s = (1 - 2Z_s)\, ds + \sqrt{\frac{\kappa Z_s(1-Z_s)}{2}}\, dB_s$$

as in [9]. Now introduce the following stopping times:

$$S = \mathrm{Inf}\{s : Z_s \in \{0,1\}\}, \qquad T = t(S) = \mathrm{Inf}\{t : R_t \in \{0,1\}\}.$$

The counterpart of Theorem 3.1 in [9] for the case $w_1 = w_2 = 0$ is obtained as Lemma 2 in the present paper; it gives the following estimate:

$$P(S > s) \asymp \exp(-\lambda(0,0)s) = \exp\left(-\frac{\kappa - 4}{2}s\right). \tag{23}$$

It remains to transfer this estimate to deterministic values of $t$. Recall that we have $2\,dt(s) = e^{2s}Z_s(1 - Z_s)\,ds$. This already proves that $dt(s) \leq e^{2s}/8\,ds$, that is, $t(s) \leq e^{2s}/16$ or $s \geq \log(16t(s))/2$. Hence,

$$P(T > t) \leq P\left(S > \frac{\log(16t)}{2}\right) \asymp t^{-(\kappa - 4)/4}.$$

To obtain the lower bound, note that the proof of Theorem 3.1 in [9] also gives the distribution of $R_s$ knowing that $S > s$, which is the eigenfunction associated to the eigenvalue $\lambda(0,0)$ for the generator of $R$, namely,

$$c[x(1 - x)]^{(\kappa - 4)/\kappa}.$$

In particular, conditionally to the fact that $S > s$, there is a positive probability that $Z_s \in [\frac{1}{4}, \frac{3}{4}]$. Comparison with Brownian motion then shows that

$$P(\forall\, s \in [s_0, s_0 + 1], Z_s \in [\tfrac{1}{8}, \tfrac{7}{8}] | Z_{s_0} \in [\tfrac{1}{4}, \tfrac{3}{4}]) \geq c > 0,$$

and combining this with (23) provides, for all $s_0 > 42$,

$$P(\forall\, s \in [s_0 - 1, s_0], Z_s \in [\tfrac{1}{8}, \tfrac{7}{8}] | S > s_0) \geq c > 0.$$

Now on this event, we obtain

$$t(s_0) \geq \int_{s_0 - 1}^{s_0} \frac{e^{2s}}{128}\,ds \geq c_0 e^{2s_0},$$

from which the lower bound follows:

$$P(T > t) \geq cP\left(S > \frac{\log(t/c_0)}{2}\right) \geq ct^{-(\kappa - 4)/4}. \qquad \square$$

The end of the proof is exactly the same as that of the previous theorem, so we do not repeat it here. $\square$

REMARK 3. For $\kappa = 8$ (where the obtained dimension is 0), the method is inconclusive; but it is possible to prove that there are no cut-points in that case; one of the simplest way to do it being to use the *reversibility* of the trace, which in turn is a consequence of the fact that $\text{SLE}_8$ is the scaling limit of the uniform spanning tree Peano path (cf. [13]). Notice however that the fact that the $\text{SLE}_8$ trace is a Peano curve is not sufficient in itself to conclude.

If $\kappa = 6$, we get that the dimension of cut-times is $1/2$. It is known in this case (using Brownian exponents) that the dimension of cut-points is $2 - 5/4 = 3/4$ (cf. [8]). For the other values of $\kappa \in (4, 8)$, the dimension of $C$ is not known.

**Acknowledgments.** I wish to thank Greg Lawler, Oded Schramm and Wendelin Werner for very useful discussions.

Departement de Mathematiques  
Batiment 425  
Universite de Paris-Sud  
F-91405  
Orsay Cedex  
France  
e-mail: vincent.beffara@math.u-psud.fr